\documentclass[10pt]{amsart}
\theoremstyle{plain}

\def\qed{{\hfill $\Box$}}

\def\Z{{\mathbb Z}}
\def\C{{\mathbb C}}

\newtheorem{thm}{Theorem}[section]
\newtheorem{cor}{Corollary}[section]
\newtheorem{prop}{Proposition}[section]

\newtheorem{lem}{Lemma}[section]
\theoremstyle{definition}
\newtheorem{defn}{Definition}[section]
\theoremstyle{remark}

\newtheorem{rem}{Remark}[section]
\newtheorem{ex}{Example}[section]

\begin{document}
\title[Representations of $U_{q}(sl_{2})$]{On Irreducible weight 
representations of a new deformation $U_{q}(sl_{2})$ of $U(sl_{2})$}
\author{Xin Tang}
\address{Department of Mathematics \& Computer Science\\
Fayetteville State University\\
Fayetteville, NC 28301}
\email{xtang@uncfsu.edu}
\keywords{ Spectrum, Deformation of $U(sl_{2})$, Irreducible
  Representations, Hyperbolic Algebras}
\thanks{}
\date{Oct 2006}
\subjclass[2000]{Primary 17A36, 17B40; Secondary 17B30}
 \begin{abstract}
Starting from a Hecke $R-$matrix, Jing and Zhang constructed 
a new deformation $U_{q}(sl_{2})$ of $U(sl_{2})$, and studied 
its finite dimensional representations in \cite{JZ}. Especically, 
this algebra is proved to be just a bialgebra, and all finite 
dimensional irreducible representations are constructed in \cite{JZ}. 
In addition, an example is given to show that not every finite
dimensional representation of this algebra is completely reducible. 
In this note, we take a step further by constructing more irreducible
representations for this algebra. We first construct points of the 
spectrum of the category of representations over this new deformation 
by using methods in noncommutative algebraic geometry. Then applied 
to the study of representations, our construction recovers all 
finite dimensional irreducible representations as constructed 
in \cite{JZ}, and yields new families of infinite dimensional 
irreducible weight representations of this new deformation 
$U_{q}(sl_{2})$. 
\end{abstract}
\maketitle

\section*{Introduction}
Spectral theory of abelian categories was first started by Gabriel 
in \cite{G}. Gabriel defined the injective spectrum of any noetherian 
grothendieck category. The injective spectrum consists of isomorphism 
classes of indecomposible injective objects in the category. Let $R$
be a commutative noetherian ring, then the injective spectrum of the 
category of all $R-$modules is isomorphic to the prime spectrum of
$R$ as affine schemes. In particular, one can reconstruct any
noetherian commutative scheme by using the injective spectrum of 
the category of quasi-coherent sheaves of modules on it. The general 
spectrum of any abelian category was later on defined by Rosenberg 
in \cite{R}. One of the advantages of this spectrum is that it is 
defined for any abelian category. Via this spectrum, one can
reconstruct any quasi-separated and quasi-compact scheme via the 
spectrum of the abelian category of quasi-coherent sheaves of modules 
over the scheme. 

Even though spectral theory is more important for the the purpose 
of noncommutative algebraic geometry, since it provides many
nontrivial examples for the study of noncommutative algebraic 
gemoetry (\cite{R}). It also has very nice applications to
representation theory. Namely, the spectrum has a natural analogue 
of the Zariski topology and its closed points are in a one to one 
correspondence with the irreducible objects in the category. In order
to study irreducible representations, one can first study the spectrum 
of the category of all representations, then single out closed ones. 

As a specific application of spectral theory to representation theory, 
points of the spectrum have been constructed for a large family of 
algebras, which are called Hyperbolic algebras in \cite{R}. And it 
is a pure luck that a lot of `small' algebras such as $U(sl_{2})$ (or its 
quantized versions) and most deformations of it are Hyperbolic algebras. 

Jing and Zhang constructed a new deformation (which is 
still denoted by $U_{q}(sl_{2})$ in \cite{JZ}) of $U(sl_{2})$ 
starting from a Hecke $R-$matrix. Among other things, they also 
studied its finite dimensional rerpesentations in detail in 
\cite{JZ}. This algebra shares a lot of properties with $U(sl_{2})$ 
and all its finite dimensional irreducible representations are 
constructed explicitly in \cite{JZ}. Note that this new deformation 
is just a bialgebra deformation of $U(sl_{2})$, and an example is 
also constructed in \cite{JZ} to show that not all of its finite 
dimensional representations are completely reducible. So the 
representation theory of this new deformation might be slightly 
different from that of the quantized enveloping algebra of ${\frak g}=sl_{2}$. And it should be an interesting problem to construct 
new families of irreducible representations for this new deformation 
for the sake of a better understanding of this algebra. 

In this note, we are interested in the problem of constructing 
irreducible representations for this algebra. We first construct 
families of points for the spectrum of the category of representations 
of this deformation. Then applied to the study of representations, our 
construction first recovers all finite dimensional irreducible 
representations of $U_{q}(sl_{2})$ as constructed in \cite{JZ}, and 
in addition yields families of infinite dimensional irreducible 
reprenestations as well. This can also be regarded as one more nice 
application of the methods in noncommutative algebraic geometry to 
representation theory. For more details about spectral theory, 
we refer the reader to \cite{R} for a complete coverage. 

The paper is organized as follows. In Section 1, we give a very 
brief review on the spectrum of any abelian category. In Section 2, 
we review the concept of Hyperbolic algebras. In Section 3, we review 
some basic facts about the new deformation $U_{q}(sl_{2})$ introduced 
by Jing and Zhang, and prove some supplementary useful Lemmas. In 
Section 4, we construct families of points of the spectrum for 
$U_{q}(sl_{2})$. Then we will use them to construct irreducible 
representations for this new deformation $U_{q}(sl_{2})$. We will 
follow the notations in \cite{JZ}, and always denote the new 
deformation by $U_{q}(sl_{2})$. The base field will always be 
fixed to be the complex field $\C$ and $q$ is not a root of unity.

\section{Basic facts about the spectrum of any abelian category}
In this section, we are going to review some basic notions and 
facts about the spectrum of any abelian category for the purpose 
of understanding the rest of this work. First, we review the
definition of the spectrum of any abelian category, then we 
explain its applications in representation theory.

Let $C_{X}$ be an abelian category. Let $M, N \in C_{X}$ be any 
two objects in $C_{X}$; We say that $M \succ N $ if and only if $N$ 
is a subquotient of the direct sum of finite copies of $M$. It is easy 
to verify that $\succ$  is a preorder. We say $M\approx N$ if and only 
if $M\succ N$ and $N\succ M$. It is obvious that $\approx$ is an
equivalence. Let $Spec(X)$ be the family of all nonzero objects $M\in
C_{X}$ such for any nonzero subobject $N$ of $M$, $N\succ M$. The 
spectrum of any abelian category is defined in \cite{R} as follows: 
\[
{\bf Spec(X)}=Spec(X)/\approx
\]
${\bf Spec(X)}$ is called the spectrum of $C_{X}$. 

${\bf Spec(X)}$ has a natural analogue of the Zariski topology. 
And its closed points are in a one to one correspondence with 
the irreducible objects of $C_{X}$. One can reconstruct any
commutative quasi-compact and semi-separated schemes from the 
spectra of the categories of quasi-coherent sheaves over them. 
The spectrum of any abelian category is more important for the 
purpose of noncommutative algebraic geometry, since it provides 
many nontrival examples of noncommutative algebraic spaces for 
the study of noncommutative algebraic geometry. 

Spectrum also has important applications in representation theory.
Let $R$ be an algebra and $C_{X}$ be the category of all $R-$modules, 
the closed points of ${\bf Spec(X)}$ are in a one to one
correspondence to irreducible $R-$modules. While ${\bf Spec(X)}$ has 
a better functorial property than the set of closed points as in the 
case of commutative algebraic geometry. So one can study the spectrum 
via the methods in noncommutative algebraic geometry, then apply to 
representation theory (\cite{R}) to recover representations. The rest 
of this paper is a typical application of spectral theory in
represenation theory. 

\section{Hyperbolic algebras $R\{\theta,\xi\}$ and points of the
  spectrum} 
Hyperbolic algebras are studied by Rosenberg in \cite{R} and  by 
Bavula under the name of Generalized Weyl algebras in \cite{B}. 
Hyperbolic algebra structure is very convenient for the construction 
of points of the spectrum. And a lot of interesting algebras such as
$U(sl_{2})$ and its quantized versions have a Hyperbolic algebra 
structure. Points of the spectrum of the category of modules over 
these algebras have been constructed in \cite{R}. In this section, we 
review some basic facts about Hyperbolic algebras and two imporatant 
construction theorems due to Rosenberg (\cite{R}).

Let $\theta$ be an automorphism of a commutative algebra $R$; and 
let $\xi$ be an element of $R$. Then we have the following definition
from \cite{R}:
\begin{defn}
We denote by $R\{\theta,\xi\}$ the corresponding $R-$algebra 
generated by $x,y$ subject to the 
following relations:
\[
xy=\xi,\quad yx=\theta^{-1}(\xi)
\]
and 
\[
xa=\theta(a)x,\quad ya=\theta^{-1}(a)y
\]
for any $a\in R$. And $R\{\theta,\xi\}$ is called a Hyperbolic algebra 
over $R$.
\end{defn}

First, we look at some basic examples of Hyperbolic algebras:
\begin{ex}
The first Weyl algebra $A_{1}$ is a Hyperbolic algebras over
$R=\C[xy]$ with $\theta(xy)=xy+1$; $U(sl_{2})$ and its quaztized 
versions are Hyperbolic algebras too. And the reader to refered to
\cite{R} for more details about Hyperbolic algebras.
\end{ex}

Let $C_{X}=C_{R\{\theta,\xi\}}$ be the category of representations 
of $R\{\theta,\xi\}$. We denote by ${\bf Spec(X)}$ the spectrum of the
category $C_{X}$. The Hyperbolic algebra structure is very convenient 
for the construction of points of the spectrum ${\bf Spec(X)}$. In
order to construct points of ${\bf Spec(X)}$, the left prime spectrum
$Spec_{l}(R\{\theta,\xi\})$ of the Hyperbolic algebra
$R\{\theta,\xi\}$ is also defined in \cite{R}, which is proved to be 
isomorphic to ${\bf Spec(X)}$. 

From now on, we will not distinguish these two spectra from each 
other. For the left prime spectrum of the Hyperbolic algebra, we 
have the following two crucial construction theorems due to Rosenberg 
from \cite{R}:
\begin{thm}
( Thm 3.2.2.in \cite{R} )
\begin{enumerate}
\item Let $P\in Spec(R)$, and the orbit of $P$ under the 
action of the automorphism $\theta$ is infinite.
\begin{enumerate}
\item If $\theta^{-1}(\xi)\in P$, and $\xi \in P$, then the left ideal 
\[
P_{1,1}\colon=P+R\{\theta,\xi\}x+R\{\theta,\xi\}y
\]
is a two-sided ideal from $Spec_{l}(R\{\theta,\xi\})$.

\item If $\theta^{-1}(\xi)\in P$, $\theta^{i}(\xi)\notin P$ for $0\leq
  i\leq n-1$, 
and $\theta^{n}(\xi)\in P$, then the left ideal 
\[
P_{1,n+1}\colon=R\{\theta,\xi\}P+R\{\theta,\xi\}x+R\{\theta,\xi\}y^{n+1}
\]
belongs to $Spec_{l}(R\{\theta,\xi\})$.

\item If $\theta^{i}(\xi)\notin P$ for $i\geq 0$ and
  $\theta^{-1}(\xi)\in P$, then
\[
P_{1,\infty}\colon=R\{\theta,\xi\}P+R\{\theta,\xi\}x
\]
belongs to $Spec_{l}(R\{\theta,\xi\})$.

\item If $\xi \in P $ and $\theta^{-i}(\xi)\notin P$ for all $i\geq
  1$, 
then the left ideal 
\[
P_{\infty,1}\colon=R\{\theta,\xi\}P+R\{\theta,\xi\}y
\]
belongs to $Spec_{l}(R\{\theta,\xi\})$.
\end{enumerate}
\item If the ideal $P$ in (b),\, (c)\, or (d) is maximal, 
then the corresponding left ideal of $Spec_{l}(R\{\theta,\xi\})$ is maximal.

\item Every left ideal $Q \in Spec_{l}(R\{\theta,\xi\})$ such that
  $\theta^{\nu}(\xi)\in Q$ for 
a $\nu \in \Z$ is equivalent to one left ideal as defined above
  uniquely from a prime ideal $P \in Spec(R)$. The latter means that 
if $P$ and $P'$ are two prime ideals of $R$ and $(\alpha,\beta)$ and
  $(\nu,\mu)$ take values $(1,\infty),(\infty,1),(\infty,\infty)$ or
  $(1,n)$, then $P_{\alpha,\beta}$ is equivalent to $P'_{\nu,\mu}$ if 
and only if $\alpha=\nu,\beta=\mu$ and $P=P'$.
\end{enumerate}
\end{thm}
\qed

\begin{thm}
(Prop 3.2.3. in \cite{R}) 
\begin{enumerate}
\item Let $P\in Spec(R)$ be a prime ideal of $R$ 
such that $\theta^{i}(\xi)\notin P$ for $i\in \Z$ and $\theta^{i}(P)-P\neq \O$ 
for $i\neq 0$, then $P_{\infty,\infty}=R\{\theta,\xi\}P\in Spec_{l}(R\{\theta,\xi\})$.
\item Moreover, if ${\bf P}$ is a left ideal of $R\{\theta,\xi\}$ such 
that ${\bf P}\cap R=P$, then ${\bf P}=P_{\infty,\infty}$. In particular,
if $P$ is a maximal ideal, then $P_{\infty,\infty}$ is a maximal 
left ideal.
\item If a prime ideal $P'\subset R$ is such that
  $P_{\infty,\infty}=P'_{\infty,\infty}$, 
then $P'=\theta^{n}(P)$ for some integer $n$. Conversely,
$\theta^{n}(P)_{\infty,\infty}=P_{\infty,\infty}$ for any $n\in \Z$. 
\end{enumerate}
\end{thm}\qed

\section{A new deformation $U_{q}(sl_{2})$ of $U(sl_{2})$}
Starting from an $R-$matrix, Jing and Zhang constructed a new 
deformaton of $U(sl_{2})$, which is still denoted by $U_{q}(sl_{2})$
in \cite{JZ}. This new deformation is a bialgebra deformation of
$U(sl_{2})$ \cite{JZ}. In this section, we first recall the definition of 
this new deformation $U_{q}(sl_{2})$. Then we verify that
$U_{q}(sl_{2})$ has a Hyperbolic algebra structure over a polynomial 
ring in two variables. Finally, we will state and verify some
supplementary useful formulas, which will be used in the next 
section. 

Let $\C$ be the field of complex numbers and $0\neq q$ be a element 
of $\C$. Let $U_{q}(sl_{2})$ be the $\C-$algebra generated by $e,f,h$ 
subject to following relations:
\[
qhe-eh=2e,\\
hf-qfh=-2f;\\
ef-qfe=h+\frac{1-q}{4}h^{2}
\]
Though this algebra is slightly different from the usual 
quanitzed enveloping algebra of $sl_{2}$, it is still denoted 
by $U_{q}(sl_{2})$ according to \cite{JZ}, and we will follow 
their notation in this paper. It is easy to see that this new 
deformation $U_{q}(sl_{2})$ shares a lot of properties with 
$U(sl_{2})$. However, this new deformation $U_{q}(sl_{2})$ is 
just a bialgebra deformation of $U(sl_{2})$ without having a 
Hopf algebra structure. The finite dimensional irreducible 
representations of this algebra were constructed in \cite{JZ}, 
and an example was constructed to show that not every finite 
dimensional representationn is completely reducible. 

In addition, it has a Casimir element which can be constructed 
as follows:  
\[
C=ef+fe+\frac{1+q}{4}h^{2}\\
=2qfe+h+\frac{1}{2}h^{2}\\
=2ef-h+\frac{q}{2}h^{2}.
\] 

We have the following basic lemma about this Casmir element $C$ 
from \cite{JZ}:
\begin{lem}
(Lem 3.4 in \cite{JZ}) $C$ $q-$commutes with generators of $U_{q}(sl_{2})$ 
in the following sense:
\begin{gather*}
eC=qCe,
\\
fC=q^{-1}Cf,
\\
hC=Ch
\end{gather*}
\end{lem}
\qed

We have the following corollary of Lemma 3.1.
\begin{cor}
$h$ commutes with $ef$.
\end{cor}
{\bf Proof:} This follows directly from the definition of $C$ and 
Lemma 3.1.
\qed

Let us denote $ef$ by $\xi$ and $e,f$ by $x,y$ respectively. 
Let $R=\C[\xi, h]$ be the subalgebra of $U_{q}(sl_{2})$ generated 
by $\xi, h$, then $R$ is a polynomial ring in two variables $\xi, h$,
which is thus commutative. We will verify that $U_{q}(sl_{2})$ is a 
Hyperbolic algebra over $R$. 

First of all, let us define an autmorphism $\theta$ of $R$ by 
\[
\theta(h)=qh-2,\quad
\theta(\xi)=q\xi+q^{2}h+\frac{q^{2}-q^{3}}{4}h^{2}-(q+1)
\]
It is obvious that $\theta$ can be extended to an algebra 
automorphism of $R$. 

In addition, we have the following basic lemma:
\begin{lem} The following identities hold:
\begin{enumerate}
\item $xh=\theta(h)x,\quad x\xi=\theta(\xi)x;$\\
\item $yh=\theta^{-1}(h)y,\quad y\xi=\theta^{-1}(\xi)y;$\\
\item $xy=\xi,\quad yx=\theta^{-1}(\xi).$
\end{enumerate}
\end{lem}
{\bf Proof:} The verification of the above lemma is straightforward,
and we will omit it. 
\qed

From Lemma 3.2., we have the following proposition:
\begin{prop}
$U_{q}(sl_{2})$ is a Hyperbolic algebra over $R$.
\end{prop}
{\bf Proof:} This follows directly from the definition of 
Hyperbolic algebras and Lemma 3.2.
\qed

\begin{cor}
The Gelfand-Kirillov dimension of $U_{q}(sl_{2})$ is $3$.
\end{cor}
{\bf Proof:} This follows from the fact that $U_{q}(sl_{2})$ is a 
Hyperbolic algebra over a polynomial algebra in two variables. Since
the Gelfand-kirillov dimension is $2$, so the Gelfand-Kirillov
dimension of $U_{q}(sl_{2})$ is $3$. 
\qed

Before we finish this section, we would like to state another useful 
lemma which will be used in the next section.
\begin{lem}
\[
\theta^{n}(h)=q^{n}h-2\frac{q^{n}-1}{q-1}
\]
and 
\[
\theta^{n}(\xi)=q^{n}\xi+\frac{q^{n+1}(1-q^{n})}{4}h^{2}+
\frac{q^{n+1}(q^{n}-1)}{q-1}h-\frac{(q^{n}-1)(q^{n+1}-1)}{(q-1)^{2}}
\]
for any $n\in \Z$.
\end{lem}
{\bf Proof:} First of all, we prove the statement is true for $n\in
\Z_{\geq 0}$. When $n=0$, the statment is obviously true. And we have 
$\theta(h)=qh-2$. Suppose the statement is true for $n-1$, we prove it 
is true for $n$ as well. Note that we have 
\[\theta^{n}(h)=\theta(\theta^{n-1}(h))=
\theta(q^{n-1}h-2\frac{q^{n-1}-1}{q-1})=q^{n}h-2\frac{q^{n}-1}{q-1}
\]
So we have proved the first statement for $n\geq 0$ by using
induction. Similar arguement shows that the statement is true for all
$n\in \Z$.

Now we are going to prove the second statement. Since $C$ is in $R$, 
then we have $xC=\theta(C)x$ and $xC=qCx$ by Lemma 3.1. So we have 
$\theta(C)=qC$. Hence $\theta^{n}(C)=q^{n}C$. Thus 
\[
2q^{n}\xi-q^{n}h+\frac{q^{n+1}}{2}h^{2}=
2\theta^{n}(\xi)-(q^{n}h-2\frac{q^{n}-1}{q-1})+
\frac{q}{2}(q^{n}h-2\frac{q^{n}-1}{q-1})^{2}
\]

Therefore, we have
\[
\theta^{n}(\xi)=q^{n}\xi+\frac{q^{n+1}(1-q^{n})}{4}h^{2}+
\frac{q^{n+1}(q^{n}-1)}{q-1}h-\frac{(q^{n}-1)(q^{n+1}-1)}{(q-1)^{2}}
\]
for any $n\in \Z$. So we are done with the proof.
\qed

\section{Construction of points of the spectrum for $U_{q}(sl_{2})$}
In this section, we construct families of points of the spectrum of
the category of representations of $U_{q}(sl_{2})$ using the
construction theorems quoted in Section 2 from \cite{R}. As a result, 
we will obtain families of irreducible weight representations of
$U_{q}(sl_{2})$. 

First of all, we have the following basic proposition:
\begin{prop}
Let $P=(\xi-\alpha,h-\beta)$ be a closed pont of $Spec(R)$. Then 
$\{\theta^{n}(p)\mid n \in \Z \}$ is a finite set if and only if 
$\alpha=-\frac{1}{(q-1)^{2}}\, \text{and} \, \beta=\frac{2}{q-1}$.
\end{prop}
{\bf Proof:} If $\alpha=-\frac{1}{(q-1)^{2}},\beta=\frac{2}{q-1}$, then
$\theta^{n}(h-\beta)=q^{n}h-2\frac{q^{n}-1}{q-1}-\beta=
q^{n}(h-\beta)+(q^{n}-1)\beta-2\frac{q^{n}-1}{q-1}$.
In addition, $\theta^{n}(\xi-\alpha)=q^{n}(\xi-\alpha)+(\frac{q^{n+1}(1-q^{n})}{4}h+
\frac{q^{n+1}(q^{n}-1)}{q-1}+\frac{q^{n+1}(1-q^{n})\beta}{4})(h-\beta)$
$+\frac{q^{n+1}(1-q^{n})\beta^{2}}{4}+\frac{q^{n+1}(q^{n}-1)\beta}{q-1}
+(q^{n}-1)\alpha-\frac{(q^{n+1}-1)(q^{n}-1)}{(q-1)^{2}}$.
So the orbit of $P$ is finite if and only if
$(q^{n}-1)\beta-2\frac{q^{n}-1}{q-1}=0$ 
and
$\frac{q^{n+1}(1-q^{n})\beta^{2}}{4}+\frac{q^{n+1}(q^{n}-1)\beta}{q-1}+
(q^{n}-1)\alpha-\frac{(q^{n+1}-1)(q^{n}-1)}{(q-1)^{2}}=0$, hence if
and only if $\alpha=-\frac{1}{(q-1)^{2}},\beta=\frac{2}{q-1}$.
\qed

For the rest of this section, we may assume that $\beta\neq \frac{2}{q-1}$.
And we have the following theorem:
\begin{thm}
Suppose that $q$ is not a root of unity and $s\in \frac{\Z_{\geq 0}}{2}$ be a
half non-negative integer. Let $P=M_{\alpha,\beta}=(\xi-\alpha,h-\beta)=(\xi-\frac{q^{-2s}-1}{1-q^{2s}},
h-2\frac{1 \mp q^{-s}}{1-q})$ be a maximal ideal of $R$, then the
corresponding point
$P_{1,n+1}=R\{\theta,\xi\}P+R\{\theta,\xi\}x+R\{\theta,\xi\}y^{n+1}$ 
of the left prime spectrum $Spec_{l}(R\{\theta, \xi\})$ a closed point. Hence 
the representation $R\{\theta,\xi \}/P_{1,n+1}$ corresponding to this 
point is a finite dimensional irreducible representation of $U_{q}(sl_{2})$.
\end{thm}
{\bf Proof:} If $P=(\xi-\alpha,h-\beta)=(\xi-\frac{q^{-2s}-1}{1-q^{2s}},h-2\frac{1\mp 
q^{-s}}{1-q})$, then we have $\theta^{-1}(\xi)\in P$ and 
$\theta^{2s}(\xi)\in P$. Thus the statement follows from (b) of part
(1) of Theorem 3.1.
\qed
\begin{rem}
The representations constructed above recover all finite dimensional
irreducible representations as constructed in \cite{JZ}.
\end{rem}

Now we are going to construct some new families of infinite
dimensional irreducible weight representations. Suppose 
$P=M_{\alpha,\beta}=(\xi-\alpha,h-\beta)$ is a maximal 
ideal of $R$, then we have the following:
\begin{thm}
\begin{enumerate}
\item If $\alpha=\beta-\frac{q-1}{4}\beta^{2}$ 
and $\beta\neq 2\frac{1\mp q^{-s}}{1-q}$ for any non-negative half integer
$s$, then the corresponding point $P_{1,\infty}\colon=R\{\theta,\xi\}P+R\{\theta,\xi\}x $ of the spectrum is closed. And the corresponding representation
$R\{\theta,\xi\}/P_{1,\infty}$ is an infinite dimensional irreducible 
highest weight representation.

\item If $\alpha=0$ and $\beta\neq 2\frac{1\mp q^{-(s+\frac{1}{2})}}{1-q}$
for any half positive integer $s$, then the corresponding point
$P_{\infty,1}\colon=R\{\theta,\xi\}P+R\{\theta,\xi\}y$ of the 
spectrum is closed, and the corresponding representation
$R\{\theta,\xi\}/P_{\infty,1}$ is an infinite dimensional 
irreducible lowest weight representation. 
\end{enumerate}
 \end{thm} 
{\bf Proof:} We will only verify the first part of this 
statement, and the rest is similar. According to Lemma 3.4, 
we have 
\[
\theta^{-1}(\xi)=q^{-1}\xi-q^{-1}h-\frac{1-q}{4q}h^{2}
\]
If $\alpha=\beta-\frac{q-1}{4}\beta^{2}$, and 
$\beta\neq 2\frac{1\mp q^{-s}}{1-q}$, then we
have $\theta^{-1}(\xi)\in P$ and $\theta^{n}(\xi)\notin P$ for any
$n\geq 0$, so that $P_{1,\infty}\colon=R\{\theta,\xi\}P+R\{\theta,\xi\}x$ 
is a closed point of the left prime spectrum
$Spec_{l}(R\{\theta,\xi\})$ by Theorem 3.1, hence the corresponding 
representation $R\{\theta,\xi\}/P{1,\infty}$ is an infinite dimensional 
highest weight irreducible representation.
\qed

\begin{thm}
Let $M_{\alpha,\beta}=(\xi-\alpha,h-\beta)$ be a maximal ideal 
of $R$ such that $\alpha\neq
\frac{q(q^{n}-1)}{4}\beta^{2}-\frac{q(q^{n}-1)}{q-1}\beta
+\frac{(q^{n}-1)(q^{n+1}-1)}{(q-1)^{2}q^{n}}$ for any $n\in \Z$, 
then the point $P_{\infty,\infty}=R\{\theta,\xi \}P\in 
Spec_{l}(R\{\theta,\xi\})$ is a closed point of the left prime 
spectrum, and the corresponding representation 
$R\{\theta, \xi\}/P_{\infty,\infty}$is an infinite dimensional 
irreducible weight representation.
\end{thm}
{\bf Proof:} The proof is a direct verification of the conditions in
Theorem 3.2, and we will omit it here.
\qed
\begin{rem}
It is tempting to construct some nonweight irreducible 
representations for $U_{q}(sl_{2})$ as done in \cite{T1} 
and \cite{T2}. Unfortunately, the Whittaker model does 
not work here. The difficulty lies in that the algebra $U_{q}(sl_{2})$
has a trivial center. So it would be an interesting problem 
to find a way of constructing nonweight irreducible 
representations for this new deformation $U_{q}(sl_{2})$.
\end{rem}

\end{document}